\documentclass[11pt,twoside,a4paper]{article}

\usepackage{mathrsfs,graphicx}
\usepackage{geometry}
\usepackage{amsmath}
\usepackage{amsfonts}
\usepackage{amsthm}
\usepackage{colortbl}

\newtheorem{lemma}{Lemma}
\newtheorem{theorem}{Theorem}
\newtheorem{remark}{Remark}
\newtheorem{proposition}{Proposition}

\newcommand{\RR}{\mathbb R}

\newcommand{\bq}{\begin{equation}}
	\newcommand{\eq}{\end{equation}}
\newcommand{\baq}{\begin{flalign}}
	\newcommand{\eaq}{\end{flalign}}
\newcommand{\bqn}{\begin{*flalign}}
	\newcommand{\eqn}{\end{*flalign}}

\newcommand{\lt}{\left}
\newcommand{\rt}{\right}
\newcommand{\lal}{\langle}
\newcommand{\ral}{\rangle}

\usepackage{authblk}
\title{Collisionless and Decentralized Formation Control for Strings}
\author[1]{Young-Pil Choi}
\author[2]{Dante Kalise}
\author[3]{Andr\'es A. Peters}

\affil[1]{Department of Mathematics, Yonsei University, Seodaemun-Gu, Seoul 03722, Republic of Korea}
\affil[2]{Department of Mathematics, Imperial College London, South Kensington Campus, SW7 2AZ London, United Kingdom}
\affil[3]{Faculty of Engineering and Sciences, Universidad Adolfo Ib\'a\~nez, Santiago 7941169, Chile}

\date{}                     
\begin{document}



\maketitle


\begin{abstract}
A decentralized feedback controller for multi-agent systems, inspired by vehicle platooning, is proposed. The closed loop resulting from the decentralized control action has three distinctive features: the generation of collision-free trajectories, flocking of the system towards a consensus state in velocity, and asymptotic convergence to a prescribed pattern of distances between agents. For each feature, a rigorous dynamical analysis is provided, yielding a characterization of the set of parameters and initial configurations where collision avoidance, flocking, and pattern formation are guaranteed. Numerical tests assess the theoretical results presented.    
\end{abstract}


\section{Introduction}
Multi-agent systems (MAS) have proven to be a versatile framework for studying diverse scalability problems in Science and Engineering, such as dynamic networks \cite{olfati}, autonomous vehicles \cite{balch}, collective behaviour of humans or animals \cite{sumpter,vicsek}, and many others \cite{surveyalbi,beaver2021overview}. Mathematically, MAS are often modelled as large-scale dynamical systems where each agent can be considered as a subset of states, updated via interaction forces such as attraction, repulsion, alignment, etc., \cite{helbing,cs07}  or through the optimization of a pay-off function in a control/game framework \cite{lasry,huang}.

In this work, we approach the study of MAS from a control viewpoint. We study a class of sparsely interconnected agents in one dimension, interacting through nonlinear couplings and a decentralized control law. The elementary building block of our approach is the celebrated Cucker-Smale model for consensus dynamics \cite{cs07}, which corresponds to a MAS where each agent is endowed with second-order nonlinear dynamics for velocity alignment, and where the influence of neighbouring agents decays with distance. The Cucker-Smale model and variants can represent the physical motion of agents on the real line, inspired by autonomous vehicle formations in platooning with a nearest-neighbour interaction scheme \cite{sepahe04,wang2019survey}. The couplings to be studied are motivated by the more general setting of the Cucker-Smale dynamics in arbitrary dimension. The original Cucker-Smale dynamics considers full network connectivity in the agent interactions, generating flocking dynamics capable of exhibiting emergent consensus behavior, that is, agents that may reach a common velocity in steady state, without the action of external forces. This framework has been extended in several directions, being most notable the inclusion of forcing terms and control \cite{cuckerhuepe,HaHaKim,cfpt,hk18}, optimal control \cite{bbck,bfk,AKKJMLR}, formation control \cite{pkh10,perea,ckpp19,zhang2023pattern}, leadership \cite{dalmao}, graph topologies \cite{ha20critical}, stochasticity \cite{huang2022stochastic}, short-range interactions \cite{yin2022nonexistence} and collision-avoidance capabilities \cite{CD,achl,bf,Choi2017,yin2020asymptotic,cheng2022collision,ha22mathphys}, the latter being an increasingly sought after property of formation control schemes for autonomous fleets of vehicles.

The main contribution of this paper is to propose a nearest neighbour interaction of agents on the real line which exhibits emergent consensus and collision avoidance under the action of a simple decentralised control law. The proposed feedback enforces a desired steady-state inter-agent distance and is inspired by formation control in vehicular platoons \cite{sepahe04,Peters}. Such a model has the potential to achieve the goals of platooning applications, that is, the coordinated, scalable, secure and efficient travel of automated vehicles \cite{feng2019}, while at the same time offering flocking, pattern formation, and collision avoidance features from the non-linear dynamics. Moreover, we provide collision-avoidance guarantees that, from a safety viewpoint, do not rely on traditional concepts in platooning such as string stability.
This work aims to provide alternative techniques for collision avoidance and formation control in platooning, that avoid the need for long-range wireless networks while considering simple and interpretable control actions.

For the derivation of collision-avoidance results, we consider the framework developed in \cite{cch14,CCMP,mp18}, which uses singular interaction kernels that blow up whenever two agents are located at the same position. We modify this setup to also consider agents with a \emph{volume} by including a threshold inter-agent distance where the kernel becomes singular.

Our main result is a rigorous characterization of flocking, collision-avoidance, and platooning behaviour for the proposed nonlinear model, in terms of the initial configuration of the system, interaction and control law parameters. Recent works have studied different aspects of the one-dimensional case \cite{piccoli2015control,ha2018first,ha2020critical,kim2021first,ZHANG2020201,choi2021one,byeon2023emergence,leslie2024finite}, relaying on full connectivity of the agent network. Our contribution differentiates itself from these works showing that these emergent properties are still present when a highly sparse nearest neighbour interaction is considered in the 1D case.

The remainder of the paper is structured as follows. In Section \ref{ps} we present the proposed model to be studied, a Cucker-Smale model with nearest neighbor singular interactions and a decentralized feedback control. We also define here a total energy functional $E(x,v)$ for the model and show that it is not increasing in time. In Section \ref{regularity} we give the results ensuring the collision-avoidance behaviour of the controlled system, and Section \ref{flocking} includes a flocking estimate showing that the velocity alignment between individuals and the inter-agent distances are uniformly bounded in time. In Section \ref{formation} we present the main formation control result. We provide different numerical experiments that illustrate our theoretical results in Section \ref{numerics} along with some concluding remarks.

\section{Problem description and preliminary results}\label{ps}

We consider a string of $N$ agents each characterized by a pair $(x_i(t),v_i(t))$ in $\RR^2$ evolving in time $t$ through second-order dynamics of the form
\begin{align} \label{main_eq}
	\left\{	
	\begin{aligned}	
		\frac{d x_i(t)}{dt}&= v_i(t),\quad i=1,\dots, N, \quad t > 0,\cr
		\frac{d v_i(t)}{dt}&=I_i(t) + u_i(t),
	\end{aligned}\right.
\end{align}
subject to initial data
\begin{flalign} \label{ini_main_eq}
	& (x_i(0), v_i(0)) =: (x_i^0, v_i^0) \quad \mbox{for} \quad i =1,\dots,N.
\end{flalign}
\begin{figure}[!ht]
	\centering
	\includegraphics[width=0.7\textwidth]{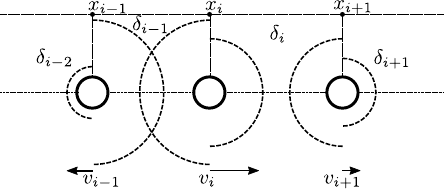}
	\caption{Diagram at a particular instant of three consecutive agents for the considered MAS. The singular interactions, providing \emph{barriers} to the agents, are indicated with the radii $\delta_i$ of the semicircles. Note that with this nomenclature, in steady state, the agents will not be considered to have collided whenever $|x_{i-1}-x_i|-\delta_{i-1}>0$. Also note that the first and last agent only have a barrier on one of their sides}
	\label{fig_setup}
\end{figure}
Here, the term $I_i$ describes nonlocal velocity interactions between individuals which are weighted by a singular communication function $\psi(r):\RR^+\longrightarrow\RR$,
\begin{align*}
	\left\{	
	\begin{aligned}	
		&I_1 =  \psi(|x_2 - x_1| - \delta_1)(v_2 - v_1),&\cr
		&I_k =  \psi(|x_k - x_{k-1}|- \delta_{k-1})(v_{k-1} - v_k) &\cr
		&\hphantom{I_N=}+ \psi(|x_{k+1} - x_k|- \delta_k)(v_{k+1} - v_k),&\cr
		&I_N = \psi(|x_N - x_{N-1}|- \delta_{N-1})(v_{N-1} -  v_N),&
	\end{aligned}\right.\qquad k=2,\dots, N-1,
\end{align*}
where the parameters $\delta_i > 0$ are fixed. This interaction term induces consensus in the velocities of the agents while preventing collisions. The second term $u_i$ serves as a decentralized feedback control depending on weight function $\phi(r):\RR^+\longrightarrow\RR$ and is given by
\begin{align*}
	\left\{	
	\begin{aligned}		
		&u_1 = -\phi(|x_1 - x_2 - z_1|^2)(x_1 - x_2 - z_1),\cr
		&u_k = \phi(|x_{k-1} - x_k - z_{k-1}|^2)(x_{k-1} - x_k - z_{k-1})\\
		& \hphantom{u_k=}-\phi(|x_k - x_{k+1} - z_k|^2)(x_k - x_{k+1} - z_k),\cr
		&u_N = \phi(|x_{N-1} - x_N - z_{N-1}|^2)(x_{N-1} - x_N - z_{N-1}),
	\end{aligned}\right.\qquad k=2,\dots, N-1\,.
\end{align*}
This feedback also depends on a vector of relative distances $z := (z_1, \dots, z_{N-1}) \in \RR^{N-1}$. The objective of this control law is to induce the formation of a string pattern characterized by $z$. The complete setting is depicted in Figure \ref{fig_setup}. For the sake of clarity, the weight functions $\psi$ and $\phi$ are chosen as
\begin{align}\label{kers}
	\psi(r) = \frac{1}{r^\alpha}, \quad\text{and}\quad \phi(r) = \frac{1}{(1 + r)^\beta}, \quad \alpha, \beta > 0\,.
\end{align}
However, the results we will state in the forthcoming sections can extended to the case in which $\phi$ is bounded and Lipschitz continuous. We will establish conditions under which the string \eqref{main_eq} converge to a consensus state with a prescribed formation while avoiding collisions between agents. For this, we begin by stating an a priori energy estimate which will be significantly used for estimating consensus emergence. We first define the total energy functional
\begin{flalign*}
	E(x,v) &:= E_1(v) + E_2(x)=  \frac{1}{4N}\sum_{i,j=1}^N |v_i - v_j|^2 + \frac 12\sum_{i=2}^N\int_0^{|x_{i-1} - x_i - z_{i-1}|^2}\phi(r)\,dr,
\end{flalign*}
and its dissipation rate
\[
D(x,v) := \sum_{i=2}^{N}\psi(|x_i - x_{i-1}|- \delta_{i-1})(v_i -v_{i-1})^2
\]

\begin{lemma}\label{lem_energy}Let $\{(x_i,v_i)\}_{i=1}^N$ be a smooth solution to the system \eqref{main_eq} on the time interval $[0,T]$. Then:
	\begin{itemize}
		\item[(i)] the mean velocity is conserved in time:
		\[
		v_c(t):=\frac1N \sum_{i=1}^N v_i(t) = v_c(0)\,.
		\]
		\item[(ii)] the total energy is not increasing in time:
		\[
		\frac{d}{dt}E(x(t),v(t)) + D(x(t),v(t)) = 0\,.
		\]
	\end{itemize}
\end{lemma}
\begin{proof}
Even though the proof is almost the same with \cite[Lemma 3.1]{ckpp19}, we provide the proof for the completeness of our work. 

(i) A straightforward computation yields
\[
\sum_{i=1}^N I_i = \sum_{i=1}^N u_i =  0.
\]
This implies
\bq\label{est_v00}
\frac{d}{dt} v_c(t) = 0, \quad \mbox{i.e.} \quad v_c(t) = v_c(0)
\eq
for $t \in [0,T]$. 

(ii) Let us first begin with the estimate for the kinetic energy:
\begin{align}\label{est_v0}
	\begin{aligned}
		\frac{1}{4N}\frac{d}{dt}\sum_{i,j=1}^N (v_i - v_j)^2 &= \frac{1}{2N}\sum_{i,j=1}^N  (v_i - v_j) \lt(\frac{dv_i}{dt} - \frac{dv_j}{dt} \rt)= \frac{1}{2N} \sum_{i,j=1}^N \lt( v_i\frac{d v_i}{dt} + v_j \frac{d v_j}{dt} \rt)\cr
		&= \sum_{i=1}^N v_i \frac{d v_i}{dt}= \sum_{i=1}^N v_i(I_i + u_i),
	\end{aligned}
\end{align}
where we used \eqref{est_v00}. Here, we use the same idea of \cite[Lemma 3.1]{ckpp19} to obtain
\begin{align}\label{est_i2}
	\begin{aligned}
		\sum_{i=1}^N v_i u_i = -\frac 12\frac{d}{dt}\sum_{i=2}^N\int_0^{|x_{i-1} - x_i - z_{i-1}|^2}\phi(r)\,dr.
	\end{aligned}
\end{align}
On the other hand, we estimate the term with $I_i$ as
\begin{align*}
\begin{split}
	\sum_{i=2}^{N-1} v_i I_i = & \sum_{i=2}^{N-1} v_i \lt(\psi(|x_i - x_{i-1}|- \delta_{i-1})(v_{i-1} - v_i) \right. \left.+ \psi(|x_i - x_{i+1}|- \delta_i)(v_{i+1} - v_i)\rt)\cr
	= & \sum_{i=2}^{N-1} v_i \lt(\psi(|x_i - x_{i-1}|- \delta_{i-1})(v_{i-1} - v_i)\rt) +\ldots\\ &\sum_{i=3}^{N} v_{i-1} \lt(\psi(|x_i - x_{i-1}|- \delta_{i-1})(v_i - v_{i-1})\rt)\cr
	= & v_2 \psi(|x_2 - x_1|- \delta_1)(v_1 - v_2)   -\ldots\\&\sum_{i=3}^{N-1}\psi(|x_i - x_{i-1}|- \delta_{i-1})(v_{i-1} - v_i)^2 + v_{N-1} \psi(|x_N - x_{N-1}|- \delta_{N-1})(v_N - v_{N-1}) \cr
	= & -v_2 I_1 - v_{N-1} I_N  -\sum_{i=3}^{N-1}\psi(|x_i - x_{i-1}|- \delta_{i-1})(v_{i-1} - v_i)^2.
 \end{split}
\end{align*}
This asserts
\begin{align}\label{est_i1}
	\begin{aligned}
		\sum_{i=1}^N v_i I_i &= v_1 I_1 + \sum_{i=2}^{N-1} v_i I_i + v_N I_N \cr
  &=  - (v_2 - v_1)I_1 -\sum_{i=3}^{N-1}\psi(|x_i - x_{i-1}|- \delta_{i-1})(v_{i-1} - v_i)^2- (v_{N-1} - v_N)I_N \cr
		&=- \sum_{i=2}^{N}\psi(|x_i - x_{i-1}|- \delta_{i-1})(v_{i-1} - v_i)^2.
	\end{aligned}
\end{align}
Combining the estimates \eqref{est_v0}, \eqref{est_i1}, and \eqref{est_i2}, we conclude the desired result.
\end{proof}

%
%
%
%
\section{Global and local existence of solutions}\label{regularity}
In this section, we show that system \eqref{main_eq}, under certain parametric and initial conditions, exhibits a non-collisional behaviour, which together with Cauchy-Lipschitz theory, subsequently provides global-in-time existence and uniqueness of smooth solutions to the system \eqref{main_eq}-\eqref{ini_main_eq}.
Inspired by \cite{ckpp19}, we present two results regarding existence of non-collisional trajectories. The first theorem requires a prescribed ordering for the initial datum $x_i^0$ and the power $\alpha$ in $\eqref{kers}$ to be $\alpha \geq 1$. The second result requires $\alpha \geq 2$, but the initial ordering assumption is removed from the string. A third result characterizes a pathological case, where a 2-agent string blows up in finite time.

\begin{theorem}\label{thm_ext0} Suppose that $\alpha \geq 1$ and the initial configuration $x_0$ satisfies $x_{i+1}^0 > x_i^0 + \delta_i $ for all $i=1,\dots,N-1$. 	Then, there exists the global unique smooth solution to the system \eqref{main_eq}-\eqref{ini_main_eq} satisfying $ x_{i+1}(t) > x_i(t) + \delta_i $ for all $i=1,\dots,N-1$ and all $t > 0$.
\end{theorem}
\begin{proof} We first notice that $\psi(x_{i+1} - x_i - \delta_i)$ is regular as long as $x_{i+1} > x_i + \delta_i$, and thus there exists a unique smooth solution to the system \eqref{main_eq}. For a fixed $T \in (0,\infty)$, let us assume that there is $t_* \in (0,T]$ where the smoothness of solutions breaks down for the first time, i.e. there is an index $\ell$ such that
\bq\label{def_t}
x_{\ell+1}(t) - x_\ell(t) > \delta_\ell \mbox{ for } t \in (0,t_*) \quad \mbox{and} \quad \lim_{t \to t_* -} x_{\ell+1}(t) - x_\ell(t) = \delta_\ell.
\eq
We denote by $[\ell]$ the set of such indices and set $i_*= \min [\ell]$. We first claim $i_* \geq 2$. If $i_* = 1$, then for $t \in (0,t_*)$ we estimate 
\[
\frac{d}{dt}\Psi(x_2 - x_1 - \delta_1) = \psi(x_2 - x_1 - \delta_1)(v_2 - v_1) = I_1 = \frac{d}{dt} (v_1 - v_c) - u_1\,,
\]
where we used \eqref{est_v00} and $\Psi$ is the primitive of $\psi$, i.e.
\[
\Psi(r) = \left\{ \begin{array}{ll}
	\displaystyle \ln(r) & \textrm{for $\alpha =1$},\\[2mm]
	\displaystyle \frac{1}{1-\alpha} r^{1-\alpha} & \textrm{for $\alpha >1$}
\end{array} \right.
\]
From this, we deduce that
\bq\label{est_claim1}
\Psi(x_2(t) - x_1(t) - \delta_1) = \Psi(x_2^0 - x_1^0 - \delta_1) = (v_1(t) - v_c(t)) - (v_1^0 - v_c(0)) - \int_0^t u_1(s)\,ds
\eq
for $t \in [0,t_*)$. On the other hand, by H\"older's inequality we find
\[
|v_1(t) - v_c(t)| = \lt|\frac1N \sum_{k=1}^N(v_1(t) - v_k(t)) \rt| \leq  \sqrt{\frac1N\sum_{k=1}^N(v_1(t) - v_k(t))^2} 
\]
and
\[
|u_1(s)| \leq \|\phi\|_{L^\infty}|x_2 - x_1 - z_1| \leq \|\phi\|_{L^\infty} \lt(|z_1| + |x_2^0 - x_1^0| + \int_0^s |v_2(\tau) - v_1(\tau)|\,d\tau \rt).
\]
These observations together with the energy estimate in Lemma \ref{lem_energy} imply that the right hand side of \eqref{est_claim1} is bounded on the time interval $(0,t_*)$, and subsequently $ t \mapsto \Psi(x_2(t) - x_1(t) - \delta_1)$ is bounded on the time interval $[0,t_*)$. This is a contradiction to \eqref{def_t} and thus the claim follows.
By the definition of $i_*$, there exists a constant $c_{i_*}>0$ such that 
\bq\label{est_x_i}
x_{i_*}(t) - x_{i_*-1}(t) - \delta_{i_* - 1} > c_{i_*}
\eq
for all $t \in (0,t_*)$. Similarly as above, we now estimate
\begin{align*}
	\frac{d}{dt}\Psi(x_{i_*+1} - x_{i_*} - \delta_{i_*}) &= \psi(x_{i_*+1} - x_{i_*} - \delta_{i_*})(v_{i_*+1} - v_{i_*}) \cr
	&= I_{i_*} + \psi(x_{i_*} - x_{i_*-1} - \delta_{i_*-1})(v_{i_*} - v_{i_*-1})\cr
	&= \frac{d}{dt} (v_{i_*} - v_c) + \psi(x_{i_*} - x_{i_*-1} - \delta_{i_*-1})(v_{i_*} - v_{i_*-1}) - u_{i_*},
\end{align*}
and thus 
\begin{align}\label{est_new}
	\begin{aligned}
		&\Psi(x_{i_*+1}(t) - x_{i_*}(t) - \delta_{i_*}) \cr
		&\quad = \Psi(x_{i_*+1}^0 - x_{i_*}^0 - \delta_{i_*}) + (v_{i_*}(t) - v_c(t)) - (v_{i_*}^0 - v_c(0))\cr
		&\qquad +\int_0^t \psi(x_{i_*}(s) - x_{i_*-1}(s) - \delta_{i_*-1})(v_{i_*}(s) - v_{i_*-1}(s))\,ds - \int_0^t u_{i_*}(s)\,ds
	\end{aligned}
\end{align}
for $t \in (0,t_*)$. Here the boundedness of the second and fourth terms can be obtained by using almost the same argument as above. We also use \eqref{est_x_i} to obtain
\[
|\psi(x_{i_*}(s) - x_{i_*-1}(s) - \delta_{i_*-1})(v_{i_*}(s) - v_{i_*-1}(s)| \leq c_{i_*}^{-\alpha} 4N E_1(v(t)) \leq c_{i_*}^{-\alpha} 4N E(x^0,v^0).
\]
Hence, the right hand side of \eqref{est_new} is bounded on the time interval $[0,t_*)$, so is the left hand side. This leads to a contradiction and thus the unique smooth solution can be actually exists up to an arbitrary finite time $T>0$. This completes the proof.
\end{proof}

We next present the second existence theorem whose proof is based on the energy estimate. For this, we first introduce a function $L^{\alpha-2}_\delta$ with $\alpha \geq 2$ given by
\[
L^{\alpha-2}_\delta(t) = \left\{ \begin{array}{ll}
	\displaystyle \sum_{i=1}^{N-1}(|x_i(t) - x_{i+1}(t)| - \delta_i)^{-{(\alpha-2)}} & \textrm{for $\alpha > 2$},\\[4mm]
	\displaystyle \sum_{i=1}^{N-1} \log (|x_i(t) - x_{i+1}(t)| - \delta_i) & \textrm{for $\alpha = 2$.}
\end{array} \right.
\]
Note that $|L^{\alpha-2}(t)| < \infty$ for $t \in [0,T]$ for some $\alpha \geq 2$ if and only if the distances between 
agents $x_i(t)$ and $x_{i+1}(t)$ are strictly greater than $\delta_i$ for all $i=1,\dots,N-1$ and $t \in [0,T]$.

\begin{theorem}\label{thm_ext} Suppose that $\alpha \geq 2$ and that the initial configuration $x_0$ satisfies
	\[
	|x_i^0 - x_{i+1}^0| >  \delta_i
	\]
	for all $i=1,\dots,N-1$.
	Then, there exists the global unique smooth solution to the system \eqref{main_eq}-\eqref{ini_main_eq} where the distances between agents satisfy $|x_i(t) - x_{i+1}(t)| > \delta_i$ for all $i=1,\dots,N-1$ and all $t > 0$.
\end{theorem}
\begin{proof} For the proof, as observed above, we will show the boundedness of the function $L^{\alpha-2}_\delta$. We first introduce the maximal life-span $T_0 = T(x^0)$ of the initial datum $x^0$:
\begin{align*}
	T_0:=& \sup\lt\{ s \in \RR_+: \exists \mbox{ solution $(x(t),v(t))$ for the system \eqref{main_eq}-\eqref{ini_main_eq} in a time-interval $[0,s)$}  \rt\}&
\end{align*}
By the assumption, it is clear $T_0 > 0$. We then claim $T_0 = \infty$. First, note that it follows from Lemma \ref{lem_energy} that
\bq\label{est_ene1}
\sum_{i=1}^{N-1}\int_0^t \frac{(v_{i+1}(s) - v_i(s))^2}{(|x_i(s) - x_{i+1}(s)|- \delta_i)^\alpha} \,ds \leq E(x^0, v^0).
\eq
Let us prove the claim above by dealing with two cases separately: $\alpha = 2$ and $\alpha > 2$. \newline

(i) $\alpha = 2$: A straightforward computation gives
\begin{align*}
\begin{split}
	\lt|\frac{d}{dt} L^0_\delta(t) \rt|  &= \lt|\frac{d}{dt} \sum_{i= 1}^{N-1}\log (|x_i(t) - x_{i+1}(t)| - \delta_i)\rt| 
	\\&= \lt|\sum_{i= 1}^{N-1} \frac{(x_i(t) - x_j(t))\cdot (v_i(t) - v_j(t))}{|x_i(t) - x_j(t)|(|x_i(t) - x_{i+1}(t)| - \delta_i)}\rt|\leq
 \sum_{i= 1}^{N-1} \frac{|v_i(t) - v_j(t)|}{|x_i(t) - x_j(t)| - \delta_i}\cr
 \end{split}
\end{align*}
for $t \in [0,T_0)$. This yields
\begin{flalign*}
	&\lt| \sum_{i= 1}^{N-1}\log (|x_i(t) - x_{i+1}(t)| - \delta_i)\rt|\leq \lt| \sum_{i= 1}^{N-1}\log (|x_i^0 - x_{i+1}^0| - \delta_i)\rt| +  \sum_{i= 1}^{N-1} \int_0^t \frac{|v_i(s) - v_{i+1}(s)|}{|x_i(s) - x_{i+1}(s)|- \delta_i}\,ds.
\end{flalign*}
On the other hand, by using the H\"older inequality and \eqref{est_ene1}, we estimate
\begin{align*}
	\sum_{i= 1}^{N-1} \int_0^t &\frac{|v_i(s) - v_{i+1}(s)|}{|x_i(s) - x_{i+1}(s)|- \delta_i}\,ds \leq \sqrt{t}\sum_{i= 1}^{N-1}  \lt(\int_0^t \frac{|v_i(s) - v_{i+1}(s)|^2}{(|x_i(s) - x_{i+1}(s)| - \delta_i)^2}\,ds\rt)^{1/2}
	\leq \sqrt{t(N-1)E(x^0,v^0)}.&
\end{align*}
Thus, we obtain
\begin{align}\label{al2}
	\lt|L^0_\delta(t) \rt|\leq \lt| L^0_\delta(0) \rt| + \sqrt{t(N-1)E(x^0,v^0)},
\end{align}
for $t \in [0,T_0)$.

(ii) $\alpha > 2$: Taking the time derivative to $L^{\alpha-2}_\delta$, we get (omitting the time arguments)
\begin{flalign*}
	\frac{d L^{\alpha-2}_\delta}{dt} &= -(\alpha-2)\sum_{i= 1}^{N-1} (|x_i - x_{i+1}|-\delta_i)^{-\alpha + 1}\frac{\lal x_i - x_{i+1}, v_i - v_{i+1}\ral}{|x_i-x_{i+1}|}\cr
	&\leq C\sum_{i= 1}^{N-1} (|x_i - x_{i+1}|-\delta_i)^{-\alpha + 1}|v_i - v_{i+1}|\cr
	&\leq C\sum_{i= 1}^{N-1}\frac{1}{(|x_i - x_{i+1}|-\delta_i)^{\alpha-2}} + C\sum_{i= 1}^{N-1}\frac{|v_i - v_{i+1}|^2}{(|x_i - x_{i+1}|-\delta_i)^\alpha}\cr
	&= CL^{\alpha-2}(t) + C\sum_{i= 1}^{N-1}\frac{|v_i - v_{i+1}|^2}{(|x_i - x_{i+1}|-\delta_i)^\alpha},
\end{flalign*}
for $t \in [0,T_0)$, where we used Young's inequality. Applying Gronwall's inequality to the above, we have
\begin{align}\label{al2+}
	\begin{aligned}
		L^{\alpha-2}_\delta(t) \leq L^{\alpha-2}_\delta(0)e^{Ct} + Ce^{Ct}\sum_{i= 1}^{N-1}\int_0^{t} \frac{|v_i(s) - v_{i+1}(s)|^2}{(|x_i(s) - x_{i+1}(s)|-\delta_i)^\alpha}\,ds\leq e^{Ct}\lt(L^{\alpha-2}_\delta(0) + CE(x^0,v^0)\rt),
	\end{aligned}
\end{align}
for $t \in [0,T_0)$, due to \eqref{est_ene1}. Since the right hand sides of \eqref{al2} and \eqref{al2+} are uniformly bounded in the time interval $[0,T_0)$, the life-span $T_0$ should be infinity, i.e., $T_0 = \infty$. This completes the proof.
\end{proof}

\begin{remark}In Theorem \ref{thm_ext0}, it is crucially used the fact that the system is posed in one dimension. However, Theorem \ref{thm_ext} can also deal with higher dimensional problems, see \cite{CCMP}.
\end{remark}

We conclude this section with a negative result characterizing a pathological configuration with 2 agents where the system blows up in finite time.
\begin{theorem} Let $\alpha \in (0,1)$ and $N=2$. Furthermore, we assume that $\delta_1, z_1$, and the initial data $\{(x_i^0, v_i^0)\}_{i=1}^2$ satisfy $ \delta_1 + z_1 \geq 0$, $x_2^0 > x_1^0 + \delta_1,$ and
	\bq\label{condi_v}
	v_1^0 - v_2^0 = \frac{2}{1-\alpha}(x_2^0 - x_1^0 - \delta_1)^{1-\alpha}.
	\eq
	Then, the smoothness of solutions to the system \eqref{main_eq}-\eqref{ini_main_eq} breaks down in finite time.
\end{theorem}
\begin{proof} For the proof, it suffices to show that there exists a finite time $t_* < \infty$ such that $x_1(t_*) + \delta_1 = x_2(t_*)$. For notational simplicity, we set $x := x_2 - x_1$ and $v:= v_2 - v_1$. Then we easily find that $x$ and $v$ satisfy
$$
\left\{	
\begin{aligned}	
	\frac{d x(t)}{dt}&= v(t),\cr
	\frac{d v(t)}{dt}&=-2(I_1(t) + u_1(t)) = -2\psi(x(t) - \delta_1) v - 2\phi(|x(t) + z_1|^2)(x(t) + z_1).
\end{aligned}\right.
$$
Note that the smooth solutions exist as long as $x(t) > \delta_1$, and this and the assumption $\delta_1 + z_1 \geq 0$ imply $x(t) + z_1 \geq 0$. Since $\phi \geq 0$, this implies that
\[
\frac{d v(t)}{dt} \leq -2\psi(x(t) - \delta_1) v = -2\frac{d}{dt} \Psi(x(t) - \delta_1).
\]
Here $\Psi$ is the primitive of $\psi$, i.e. 
\[
\Psi(r) = \frac{1}{1 - \alpha} r^{1-\alpha}.
\]
We then solve the above differential inequality to get
\bq\label{eqx_diff}
\frac{d(x(t) - \delta_1)}{dt} = v(t) \leq -2\Psi(x(t) - \delta_1) = -\frac{2}{1-\alpha}(x(t) - \delta_1) ^{1-\alpha}
\eq
due to \eqref{condi_v}. We notice that the above differential inequality is sub-linear, and thus there exists $t_* < \infty$ such that $x(t_*) - \delta_1 = 0$. Indeed, we obtain from \eqref{eqx_diff} that
\[
(x(t) - \delta_1)^\alpha \leq (x^0 - \delta_1)^\alpha - \frac{2\alpha}{1-\alpha}t.
\]
Hence we have
\[
t_* \leq \frac{(x^0 - \delta_1)^\alpha}{2\alpha} (1-\alpha)\,,
\]
thus completing the proof.
\end{proof}

%
%
%
%
\section{Time-asymptotic behavior}\label{flocking}
Having characterized the well-posedness of the trajectories of the system \eqref{main_eq}, we now turn our attention to the study of flocking emergence within the controlled string. In a flocking configuration, all agents travel with the same constant velocity, and as a direct consequence the distance between agents remain constant. Let us recall our energy functionals 
\[
E_1(v)=\frac{1}{4N}\sum_{i,j=1}^N |v_i - v_j|^2  \quad \mbox{and} \quad E_2(x)=  \frac 12\sum_{i=2}^N\int_0^{|x_{i-1} - x_i - z_{i-1}|^2}\phi(r)\,dr.
\]

We provide a rigorous asymptotic flocking estimate for the system \eqref{main_eq}.

\begin{theorem}\label{thm_1} 
	Suppose that either assumptions of Theorems \ref{thm_ext0} or \ref{thm_ext} hold. Furthermore, we assume 
	\bq\label{main_as}
	\int_0^\infty \phi(r)\,dr > \frac{1}{2N}\sum_{i,j=1}^N|v_i^0 - v_j^0|^2 + \sum_{i=2}^N \int^{|x_{i-1}^0 - x_{i}^0 - z_{i-1}|^2}_0 \phi(r)\,dr .
	\eq
	Then, the string converges asymptotically towards a flocking state, that is
	\[
	\sup_{0 \leq t \leq \infty}\max_{i,j=1,\dots,N} |x_i(t) - x_j(t)| < \infty \qquad\text{and}\qquad
	\max_{i,j=1,\dots,N} |v_i(t) - v_j(t)|\to 0
	\]
	as $t\to\infty$.
\end{theorem}
\begin{proof}  From Theorems \ref{thm_ext0} or \ref{thm_ext}, it follows existence and uniqueness of a smooth solution globally in time.

{\it (Uniform-in-time boundedness):} It follows from the energy estimate in Lemma \ref{lem_energy} that
\[
E_2(x(t)) \leq E(x^0,v^0),
\]
i.e.
\bq\label{est_000}
\sum_{i=2}^N\int_{|x_{i-1}^0 - x_i^0 - z_{i-1}|^2}^{|x_{i-1}(t) - x_i(t) - z_{i-1}|^2}\phi(r)\,dr \leq \frac{1}{2N}\sum_{i,j=1}^N|v_i^0 - v_j^0|^2 \quad \mbox{for} \quad t \geq 0.
\eq
On the other hand, under our main assumptions, we can find some constant $\rho > 0$ such that 
\[
\frac{1}{2N}\sum_{i,j=1}^N|v_i^0 - v_j^0|^2 + \sum_{i=2}^N \int^{|x_{i-1}^0 - x_{i}^0 - z_{i-1}|^2}_0 \phi(r)\,dr  \leq \int_0^{\rho^2} \phi(r)\,dr. 
\]
This together with \eqref{est_000} yields
\[
0 \leq \frac{1}{2N}\sum_{i,j=1}^N|v_i^0 - v_j^0|^2 \leq \int_{|x_{i-1}^0 - x_i^0 - z_{i-1}|^2}^{\rho^2}\phi(r)\,dr
\]
for all $i=2,\dots,N$. This implies that
\bq\label{est_diff}
|x_{i-1}(t) - x_i(t) - z_{i-1}| \leq \rho \quad \mbox{for} \quad i=2,\dots,N.
\eq
For any $i < j,$ by telescoping and the triangle inequality we estimate
\begin{align*}
	|x_i - x_j| = \lt| \sum_{\ell = i}^{j-1}(x_\ell - x_{\ell+1})\rt| \leq \sum_{\ell = i}^{j-1}|x_\ell - x_{\ell+1}| \leq \sum_{\ell = i}^{j-1}|x_\ell - x_{\ell+1} - z_{\ell}| + \sum_{\ell = i}^{j-1}|z_{\ell}|,
\end{align*}	
and thus
\[
|x_i - x_j| \leq |j-i|\, \rho + \sum_{\ell = i}^{j-1}|z_{\ell}| \leq  (N-1)\, \rho + \sum_{i= 1}^{N-1}|z_i| < \infty\,,
\]
given the boundedness of distances between agents at all times.

{\it (Velocity alignment behavior):} From the bound above, we find
\[
|x_i - x_{i-1}|- \delta_{i-1} \leq |x_i - x_{i-1} - z_{i-1}| + |z_{i-1}| - \delta_{i-1} \leq \rho + |z_{i-1}| + \delta_{i-1} \leq \rho + \max_{i=1,\dots,N-1}(|z_i| + \delta_i).
\]
Since $\psi$ is monotonically decreasing, we obtain
\[
\psi_m := \min_{2 \leq i \leq N} \psi(|x_i - x_{i-1}|- \delta_{i-1}) \geq \psi\lt(\rho + \max_{i=1,\dots,N-1}(|z_i| + \delta_i)\rt) > 0.
\]
This implies that the dissipation rate $D$ is bounded from below by
\[
D(x(t),v(t)) = \sum_{i=2}^{N}\psi(|x_i - x_{i-1}|- \delta_{i-1})(v_{i-1} - v_i)^2 \geq \psi_m \sum_{i=2}^{N}(v_{i-1} - v_i)^2.
\]
Then, by Lemma \ref{lem_energy}, we get
\[
\sum_{i=2}^{N}\int_0^\infty(v_{i-1}(t) - v_i(t))^2\,dt < \infty,
\]
and subsequently, this leads to
\[
\int_0^\infty E_1(v(t))\,dt = \frac{1}{4N}\sum_{i,j=1}^N\int_0^\infty  |v_i(t) - v_j(t)|^2\,dt < \infty.
\]
Indeed, by telescoping, for any $i < j$
\[
|v_i - v_j| \leq \sum_{\ell = i}^{j-1} | v_\ell - v_{\ell+1}| \leq \sqrt{|i-j|} \sqrt{\sum_{\ell = i}^{j-1} | v_\ell - v_{\ell+1}|^2},
\]
and thus
\bq\label{est_vv}
\sum_{i,j=1}^N|v_i - v_j|^2 \leq   c_N\sum_{i=2}^N |v_{i-1} - v_i|^2\,,\quad\text{ where }\quad c_N := \sum_{i,j=1}^N |i-j|\,.
\eq
Moreover, we also find
\begin{align}\label{est_02}
	\begin{aligned}
		\lt|\sum_{i=1}^N v_i u_i \rt| &\leq \left|\sum_{i=2}^N\phi(|x_{i-1} - x_i - z_{i-1}|^2)\lal x_{i-1} - x_i - z_{i-1}, v_{i-1} - v_i \ral\right|\cr
		&\leq \rho \sum_{i=2}^N |v_{i-1}-v_i|\leq C\sqrt{E(x^0,v^0)},\cr
	\end{aligned}
\end{align}
where $C$ is independent of $t$ and we used
\[
\max_{1 \leq i,j \leq N}|v_i(t) - v_j(t)| \leq \sqrt{ \sum_{i,j=1}^N |v_i(t) - v_j(t)|^2} \leq 2\sqrt{NE(x^0,v^0)}.
\]
Furthermore, note that
\begin{flalign*}
	\frac{1}{4N}\sum_{i,j=1}^N|v_i-v_j|^2 =& \int_0^t(-D(x(s),v(s)))\,ds+ \sum_{i=1}^N\int_0^t  v_i(s) u_i(s) \,ds + \frac{1}{4N}\sum_{i,j=1}^N|v_i^0 -v_j^0|^2.
\end{flalign*}
The dissipation rate $D$ is integrable and thus the first term on the right side of the above equality is absolutely continuous. Regarding the second term, its time-derivative is uniformly bounded in time, see \eqref{est_02}, from where it follows that it is Lipschitz continuous. This implies that the $E_1(v(t))$ is the sum of an absolutely continuous function and a Lipschitz continuous function. Thus, we obtain that $E_1(v(t))$ is uniformly continuous. Since $E_1(v(t))$ is also integrable, $E_1(v(t)) \to 0$ as $t \to \infty$. This completes the proof.
\end{proof}

\begin{remark} If $\beta \leq 1$, then $\phi$ is not integrable, thus the left hand side of \eqref{main_as} becomes infinity. This implies that the assumption \eqref{main_as} automatically holds. On the other hand, if $\beta > 1$, we obtain
	\[
	\int_0^\infty \phi(r)\,dr = \int_0^\infty \frac{1}{(1 + r)^\beta}\,dr=\frac{1}{\beta-1},
	\]
	and thus \eqref{main_as} can be rewritten as 
	\begin{align}\label{eq:condT2}
		E(x_0, v_0)  < \frac{1}{2(\beta-1)}.
	\end{align}
	
	From these equivalences, it becomes evident that the fulfilment of the flocking condition depends only on two parameters of the model, namely, the number of agents $N$ in the string  and the control interaction constant $\beta>1$, which regulates the strength of the control action. The constant $\alpha$ does not play a role on the condition. Having fixed a number of agents and $\beta$, flocking solely depends on the cohesiveness of the initial configuration.
\end{remark}

\begin{remark} If we define $\Phi$ by the primitive of $\phi$, then it is clear that $r \mapsto \Phi(r)$ is strictly increasing, and thus the constant $\rho$ appeared in \eqref{est_diff} can be expressed by
	\bq\label{def_rho}
	\rho = \sqrt{ \Phi^{-1}\lt(\frac{1}{2N}\sum_{i,j=1}^N|v_i^0 - v_j^0|^2 + \Phi(|x_{i-1}^0 - x_i^0 - z_{i-1}|^2) \rt)}.
	\eq
\end{remark}

%
%
%
%
\section{Exponential emergence of pattern formation and velocity alignment}\label{formation}
In this section, we conclude our characterization of the string trajectories by studying the exponential emergence of pattern formation and velocity alignment behavior under additional assumption on the solutions. We first provide an auxiliary result, a modification of Young's inequality, which can be proved by a similar argument as in \cite[Lemma 6.1]{ckpp19}. We thus omit its proof here. 
\begin{lemma}\label{tedious}
	Let $a_1,\dots,a_{N-1}$ be a set of vectors in $\RR^d$ and $b_1,\dots,b_{N-1}$ be a set of positive scalars. Then
	\begin{align*}
		-\sum_{i=1}^{N-1}b_i |a_i|^2 + \sum_{i=1}^{N-2}b_i\lal a_i, a_{i+1}\ral \leq -\epsilon_0 \sum_{i=1}^{N-1}b_i|a_i|^2,
	\end{align*}
	where $\epsilon_0\in(0,1)$ is a sufficiently small number.
\end{lemma}

We now state our main result, which provides non-collisional behavior, flocking, and an exponential decay estimate towards the string configuration encoded in the relative distance vector $z$.

\begin{proposition}\label{pat}
	Suppose that the assumptions of Theorem \ref{thm_1} are satisfied. Furthermore, we assume that
	\bq\label{main_as2}
	\inf_{t \geq 0} \min_{1 \leq i \leq N-1} \lt(|x_i(t) - x_{i+1}(t)| - \delta_i\rt) > 0\,.
	\eq
	Then, we have 
	\[
	\max_{i=2,\dots,N} |x_{i-1}(t) - x_i(t) - z_{i-1}| + \max_{i,j=1,\dots,N} | v_i(t) - v_j(t)| \to 0
	\]
	exponentially fast as $t \to \infty$.
\end{proposition}
\begin{proof} We first notice that the energy estimate in Lemma \ref{lem_energy} only provides a dissipation rate for the velocity. In order to have a complete exponential decay estimate, it is required to obtain the dissipation rate associated to the positions. For this, we consider the following quantity:
\[
\sum_{i=1}^{N-1} ( x_i-x_{i+1}-z_i)( v_i-v_{i+1}).
\]
Note that the total energy is bounded from below and above by
\[
\frac{1}{4N}\sum_{i,j=1}^N |v_i - v_j|^2 + \phi_m \sum_{i=1}^{N-1} | x_i-x_{i+1}-z_i|^2  \leq E(x,v)  \leq \frac{1}{4N}\sum_{i,j=1}^N |v_i - v_j|^2  +  \sum_{i=1}^{N-1} | x_i-x_{i+1}-z_i|^2,
\]
where we used 
\[
\phi_m := \min_{s \in [0, \,\rho]}\phi(s) \leq \phi(r) \leq 1. 
\]
This shows that a modified total energy $E_\gamma$, defined as 
\[
E_\gamma (x,v) : = \gamma E(x,v) + \sum_{i=1}^{N-1} ( x_i-x_{i+1}-z_i)( v_i-v_{i+1}),
\]
has similar upper and lower bound estimates as the one for $E(x,v)$ when $\gamma > 0$ large enough. Indeed, for $\gamma > \sqrt{2N/\phi_m}$, we readily find
\begin{align}\label{est_eg}
	\begin{aligned}
		E_\gamma (x,v) &\geq  \frac{\gamma \phi_m}2 \sum_{i=1}^{N-1} | x_i-x_{i+1}-z_i|^2 +   \lt(\frac\gamma{4N} - \frac1{2\gamma \phi_m}\rt) \sum_{i,j=1}^N |v_i - v_j|^2\cr
		&\geq c_\gamma \lt(\sum_{i=1}^{N-1} | x_i-x_{i+1}-z_i|^2 +  \sum_{i,j=1}^N |v_i - v_j|^2\rt),
	\end{aligned}
\end{align}
where $c_\gamma > 0$ is given by
\[
c_\gamma:= \min\lt\{\frac{\gamma \phi_m}2, \frac\gamma{4N} - \frac1{2\gamma \phi_m} \rt\}.
\]
The upper bound on $E_\gamma$ can be easily obtained. On the other hand, it follows from \eqref{main_eq} that
\begin{flalign}\label{est_new0}
	\begin{aligned}
		\frac{d}{dt}\sum_{i=1}^{N-1}( x_i-x_{i+1}-z_i)( v_i-v_{i+1}) = & \sum_{i=1}^{N-1} (v_i-v_{i+1})^2 + \sum_{i=1}^{N-1}  (x_i-x_{i+1}-z_i)( I_i - I_{i+1})+ \\ &
  \sum_{i=1}^{N-1}  (x_i-x_{i+1}-z_i)( u_i -u_{i+1}).
	\end{aligned}
\end{flalign}
Since
\[
|I_1| \leq \psi_M | v_2 - v_1|, \quad |I_N| \leq \psi_M|v_{N-1} - v_N|, \quad \mbox{and} \quad |I_i| \leq \psi_M(|v_{i-1} - v_i| + |v_i - v_{i+1}|)
\]
for $i=2,\dots,N-1$, we easily find
\[
\sum_{i=1}^{N-1}|I_i - I_{i+1}|^2 \leq 2\lt(|I_1|^2 + 2 \sum_{i=2}^{N-1}|I_i|^2 + |I_N|^2 \rt) \leq 16\,\psi_M^2 \sum_{i=1}^{N-1} |v_i - v_{i+1}|^2,
\]
where $\psi_M := \sup_{r \in [0,\infty)} \psi(r) < \infty$, which can be defined by the assumption \eqref{main_as2}. Thus, by using Young's inequality we have
\bq\label{est_new1}
\sum_{i=1}^{N-1}  (x_i-x_{i+1}-z_i)( I_i - I_{i+1}) \leq \epsilon_1  \sum_{i=1}^{N-1}  |x_i-x_{i+1}-z_i|^2 + \frac{4 \psi_M^2}{\epsilon_1}\sum_{i=1}^{N-1} |v_i - v_{i+1}|^2\,,
\eq
where $\epsilon_1$ will be determined later. Regarding the term with $u_i$, we estimate
\begin{align*}
	\sum_{i=1}^{N-1}  (x_i-x_{i+1}-z_i)( u_i -u_{i+1})&= -2\sum_{i=1}^{N-1} \phi(|x_i - x_{i+1} - z_i|^2)|x_i - x_{i+1} - z_i|^2\cr
	&\quad + 2\sum_{i=1}^{N-2} \phi(|x_i - x_{i+1} - z_i|^2)(x_i - x_{i+1} - z_i)(x_{i+1} - x_{i+2} - z_{i+1}).
\end{align*}
We then use Lemma \ref{tedious} with $b_i = 2\phi(|x_i - x_{i+1} - z_i|^2)$ and $a_i = x_i - x_{i+1} - z_i$ to get
\begin{align*}
	\sum_{i=1}^{N-1}  (x_i-x_{i+1}-z_i)( u_i -u_{i+1}) &\leq -2\epsilon_0\sum_{i=1}^{N-1} \phi(|x_i - x_{i+1} - z_i|^2)|x_i - x_{i+1} - z_i|^2\cr
	&\leq -2\epsilon_0\phi_m\sum_{i=1}^{N-1}  |x_i - x_{i+1} - z_i|^2,
\end{align*}
where $\epsilon_0$ is given in Lemma \ref{tedious}. This together with \eqref{est_new0}, \eqref{est_new1}, and choosing $\epsilon_1 = \epsilon_0 \phi_m$ implies
\begin{align*}
	\frac{d}{dt}\sum_{i=1}^{N-1}( x_i-x_{i+1}-z_i)( v_i-v_{i+1}) &\leq -\epsilon_0\phi_m\sum_{i=1}^{N-1}  |x_i - x_{i+1} - z_i|^2 + \frac{4 \psi_M^2}{\epsilon_0 \phi_m}\sum_{i=1}^{N-1} |v_i - v_{i+1}|^2.
\end{align*}
Thus the modified total energy $E_\gamma$ satisfies
\[
\frac{d}{dt} E_\gamma(x,v) \leq - \lt(\gamma \psi_m - \frac{4 \psi_M^2}{\epsilon_0 \phi_m} \rt) \sum_{i=1}^{N-1}|v_i - v_{i+1}|^2 -\epsilon_0\phi_m\sum_{i=1}^{N-1}  |x_i - x_{i+1} - z_i|^2.
\]
By taking $\gamma > \max\{(4 \psi_M^2)/(\epsilon_0 \phi_m \psi_m), \sqrt{2N/\phi_m}\}$ and using \eqref{est_vv} and \eqref{est_eg}, we further estimate
\begin{align*}
	\frac{d}{dt} E_\gamma(x,v) &\leq - \frac1{c_N}\lt(\gamma \psi_m - \frac{4 \psi_M^2}{\epsilon_0 \phi_m} \rt) \sum_{i,j=1}^N|v_i - v_j|^2 -\epsilon_0\phi_m\sum_{i=1}^{N-1}  |x_i - x_{i+1} - z_i|^2\cr
	&\leq -\frac1{c_\gamma}\min\lt\{\frac1{c_N}\lt(\gamma \psi_m - \frac{4 \psi_M^2}{\epsilon_0 \phi_m} \rt), \epsilon_0\phi_m \rt\} E_\gamma(x,v).
\end{align*}
Applying Gr\"onwall's lemma to the above gives the exponential decay of the modified total energy $E_\gamma$. Moreover, the relation \eqref{est_eg} concludes the desired result.
\end{proof}
\begin{remark} The a priori assumption \eqref{main_as2} imposes some constraints on $z_i$. For instance, if we fix the order for the initial positions as $x_{i}^0 + \delta_i < x_{i+1}^0$ for $i=1,\dots,N-1$, then by Theorem \ref{thm_ext}, we have $x_i(t) + \delta_i < x_{i+1}(t)$ for all $t \geq 0$. This implies that in order to have the time-asymptotic pattern formation, $z_i$ and $\delta_i$ should satisfy $z_i > \delta_i$ for all $i=1,\dots,N-1$.
	
	On the other hand, if we assume 
	\[
	z_{i-1} > \rho + \delta_{i-1}, \quad i=2,\dots,N,
	\]
	where $\rho$ is given as in \eqref{def_rho}, then
	\[
	|x_i(t) - x_{i-1}(t)|- \delta_{i-1} \geq z_{i-1} - \rho - \delta_{i-1} > 0
	\]
	for all $t \geq 0$ and all $i=2,\dots,N$. Indeed, it follows from \eqref{est_diff} that
	\[
	|x_i(t) - x_{i-1}(t)| = |x_i(t) - x_{i-1}(t) - z_{i-1} + z_{i_1}| \geq z_{i-1} - \rho,
	\]
	thus subtracting $\delta_{i-1}$ from the both sides gives the desired result.
\end{remark}

\begin{figure}[!t]\begin{center}
		\includegraphics[width=\columnwidth]{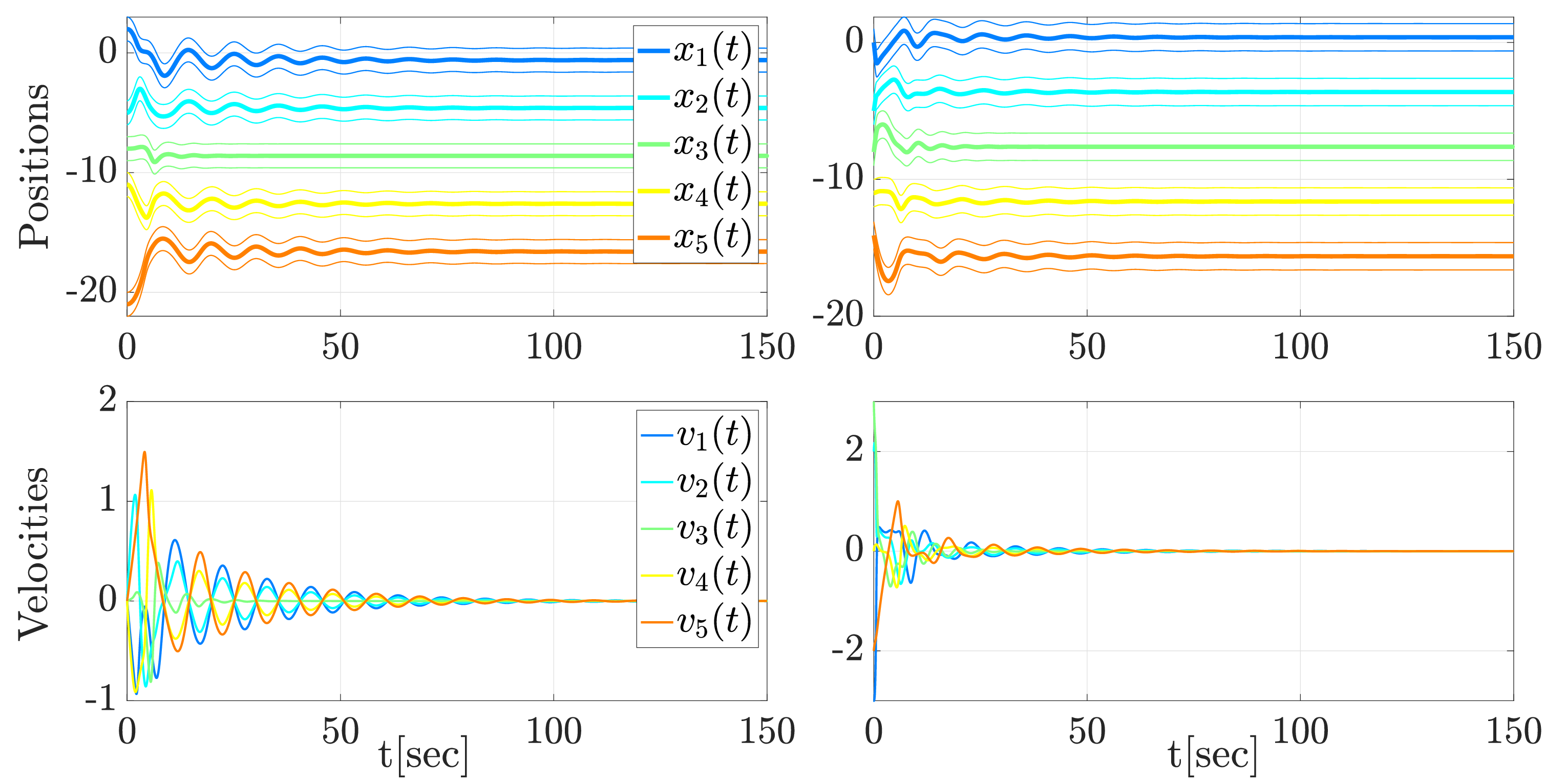}\\ 
		\caption{Positions over time of 5 particles on the line. Left: case 1) agents initially at rest and located at non-collided positions. Right: case 2) agents with non-zero initial velocities and located at close proximity but not collided. Thin lines represent the \emph{volumes} of the agents. No collisions occur due to the singular interactions and the desired formation is acquired in steady state. It can also be noted that the velocities seem to converge exponentially fast. Note that the \emph{volumes} seem to change in time but this is only due to the chosen visualization style.} \label{fig:ex1}
	\end{center}
\end{figure}

\section{Numerical examples}\label{numerics}
In the following, we further elucidate the applicability of our results through two numerical experiments illustrating string flocking to a pattern formation, as well as energy evolution. 

\subsection{Regular collision-less behavior of the interconnected system}

For simplicity in the visualization, we will first consider a collection of $N=5$ agents. We will assume that the agents are in the desired ordering, that is, the final configuration of the agents does not require a collision to occur. For the model parameters, we select $\alpha=2.1$, $\beta=0.8$ and $\delta=2$, equal for all inter-agent barriers. We will consider two cases for the agents' initial conditions: 1) the agents are initially at rest and located at non-collided positions, 2) the agents have different initial velocities with $v_c(0)=0$ and are close to each other but not collided. In both cases, the desired inter-agent spacings are given by $z_i=\delta+2$ (note that $\delta$ is added to avoid configurations that are collided in steady state). In such a scenario, the agents should reposition themselves reaching the consensus velocity of zero and the desired pattern. Figure \ref{fig:ex1} illustrates both cases. As predicted by Theorems \ref{thm_ext} and \ref{thm_1}, and given that $\alpha>2$ and $\beta< 1$ with no initial collisions, we have that the agents do not collide as time progresses, even when they are initially almost touching. Moreover, the agents reach the desired inter-agent spacings and they reach flocking exponentially fast. 

\begin{figure}[!t]\begin{center}
		\includegraphics[width=\columnwidth]{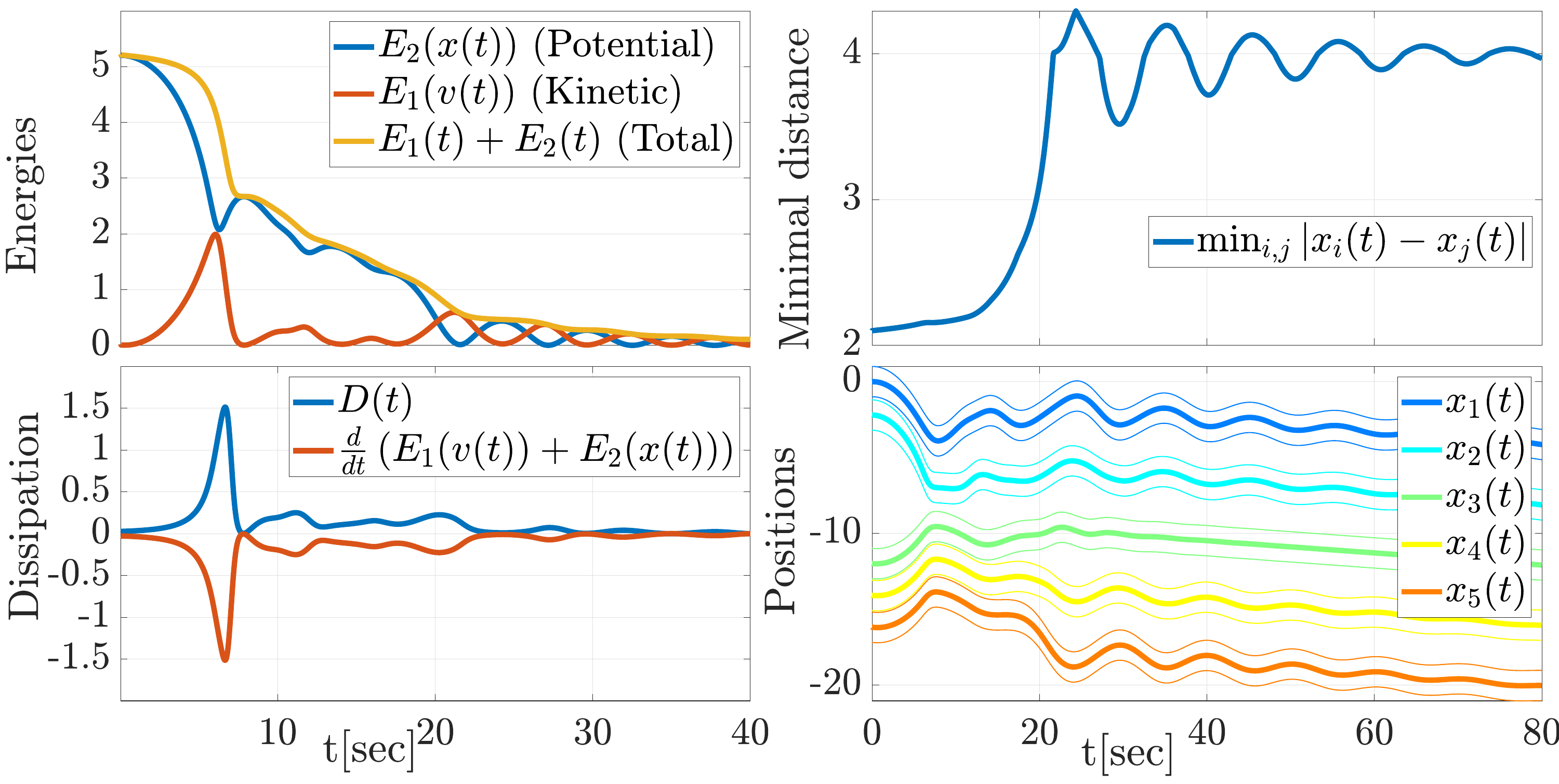}\\
		\caption{5 particles on the line with non-zero average velocity. Left-Top: Energy decomposition; Left-Bottom: Dissipation; Right-Top: Minimal distance between agents; Right-Bottom: Positions over time of the particles achieving a consensus speed and the desired spatial formation as the system evolves. Note that the energies satisfy statement (ii) in Lemma \ref{lem_energy}.} \label{fig:example_energies}
	\end{center}
\end{figure}

Now, for the same parameters as before, Figure \ref{fig:example_energies} illustrates the behaviour of the multi-agent system when $v_c(0)=-0.2$. We can observe that the statement of Lemma \ref{lem_energy} is satisfied over the time evolution of the dynamics, that is, the total energy of the system is non-increasing and the dissipation rate is entirely determined by the interaction term (see Figure \ref{fig:example_energies} bottom-left). Some agents are very close to the interaction limit determined by $\delta=2$, at $t=0[sec]$, which is highlighted by Figure \ref{fig:example_energies} top-right. The agents 2 and 3 then approach each other for a short period of time and around $t=10[sec]$ all the agents begin to spread out to reach the desired formation, with an average velocity of $-0.2$ in steady state.

\subsection{Varying barriers and small violation of the flocking conditions}

We now select the system parameters in order to study what occurs when the flocking condition in Theorem \ref{thm_1} is violated. For $N=5$, $\alpha=2.1$, $\beta=1.1$ and non-symmetric barriers given by $\delta_i=\{2.55, 2.15, 0.65, 1.05\}$. We consider the case where agents have non-zero initial velocities with $v_c(0)=0$ and are initially positioned at non-colliding locations. The desired inter-agent spacings are given by $z_i=\{6,6,6,3\}$ 

According to \eqref{eq:condT2}, since $\beta > 1$, flocking is guaranteed when $E(x_0, v_0) < \frac{1}{2(\beta-1)}$. For our chosen parameters, the theoretical threshold is $\frac{1}{2(1.1-1)} = 5$, but our initial configuration yields $E(x_0, v_0) \approx 5.2$, slightly exceeding this bound.

Despite violating the energy condition, the numerical results show that agents still achieve velocity consensus and converge to the desired formation. This does not contradict Theorem \ref{thm_1}, as the theorem provides a sufficient but not necessary condition for flocking. Collisions are still avoided through the singular interactions while the system successfully reaches the desired pattern formation in steady state.

\begin{figure}[!ht]\begin{center}
		\includegraphics[width=\columnwidth]{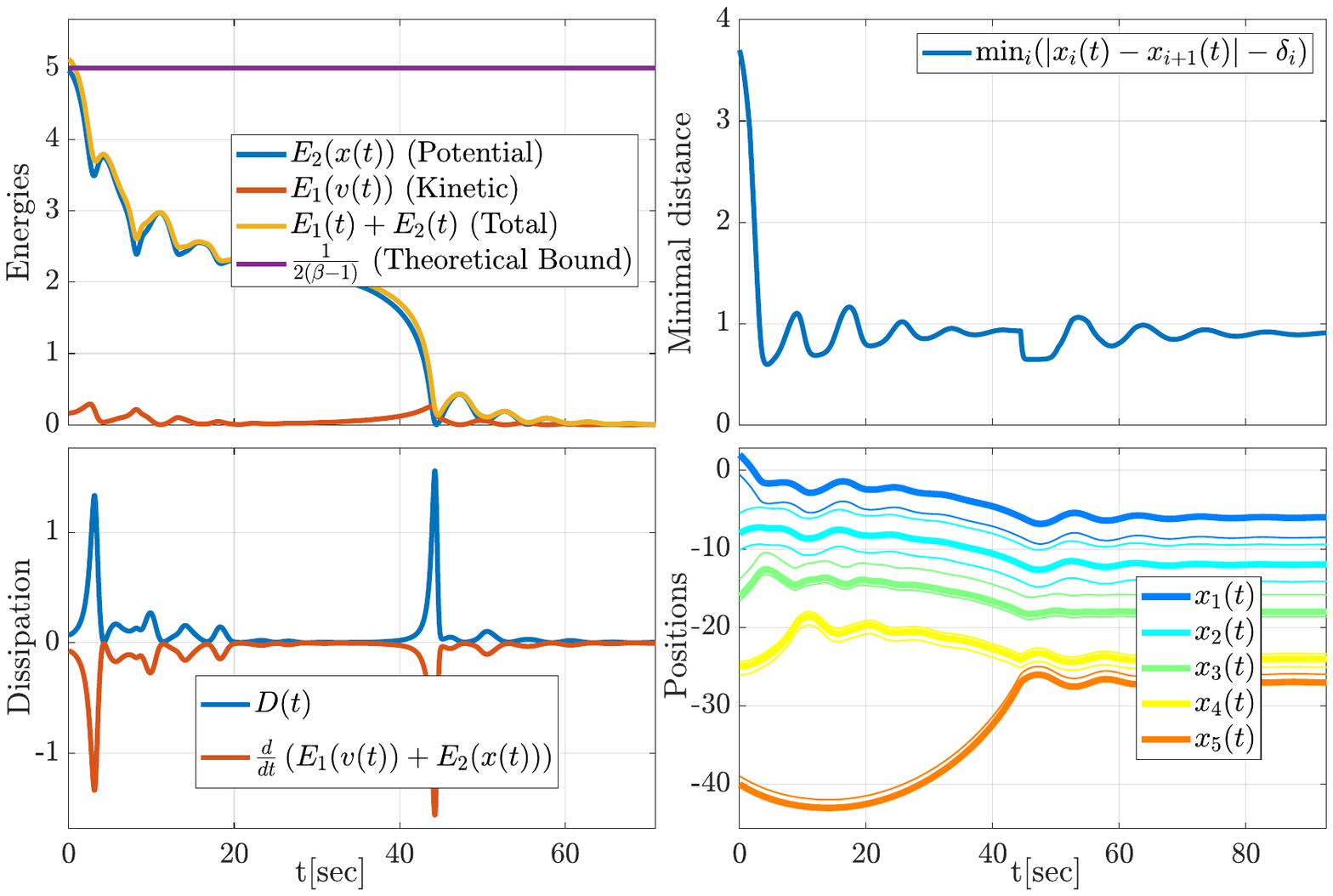}\\
		\caption{5 particles on the line with zero average velocity. Left-Top: Energy decomposition; Left-Bottom: Dissipation; Right-Top: Minimal distance between agents; Right-Bottom: Positions over time of the particles achieving a consensus speed of zero and the desired spatial formation as the system evolves. Note that the initial configuration violates \eqref{eq:condT2}.} \label{fig:example_barriers}
	\end{center}
\end{figure}

\subsection{Flocking and formation acquisition for large violation of the flocking condition}

\begin{figure}[!t]\begin{center}
		\includegraphics[width=\columnwidth]{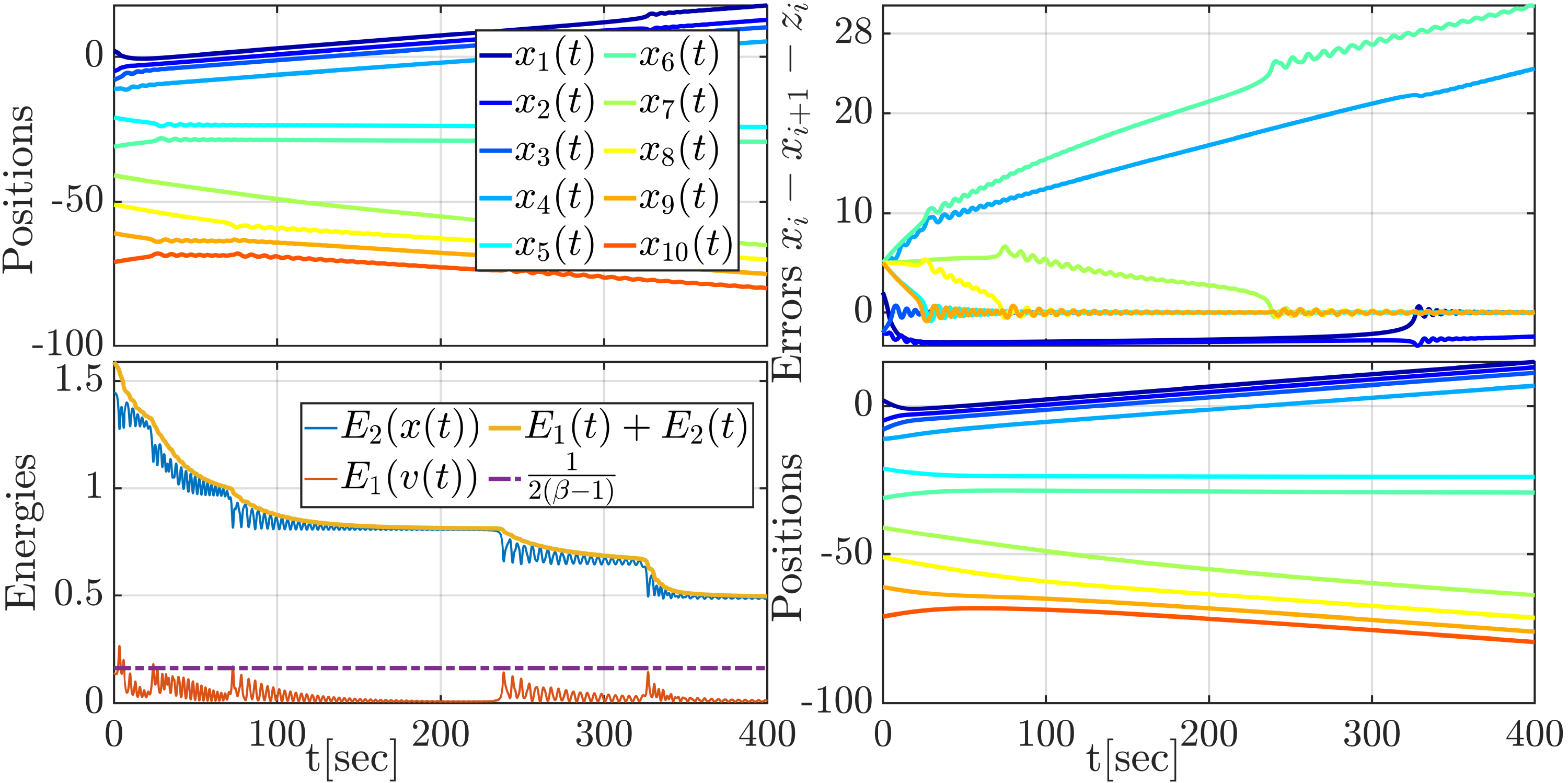}\\
		\caption{10 particles on the line with zero average velocity when $\beta=4.1$. Left-Top: Positions over time of the particles when the control is used; Left-Bottom: Energy decomposition and flocking condition for $\beta>1$; Right-Top: Errors from the desired formation $x_i-x_{i+1}-z_i$; Right-Bottom: Positions over time of the particles when the control is not used. Flocking does not occur, although collisions are still avoided.} \label{fig:beta4.1}
	\end{center}
\end{figure}

\begin{figure}[!t]\begin{center}
		\includegraphics[width=\columnwidth]{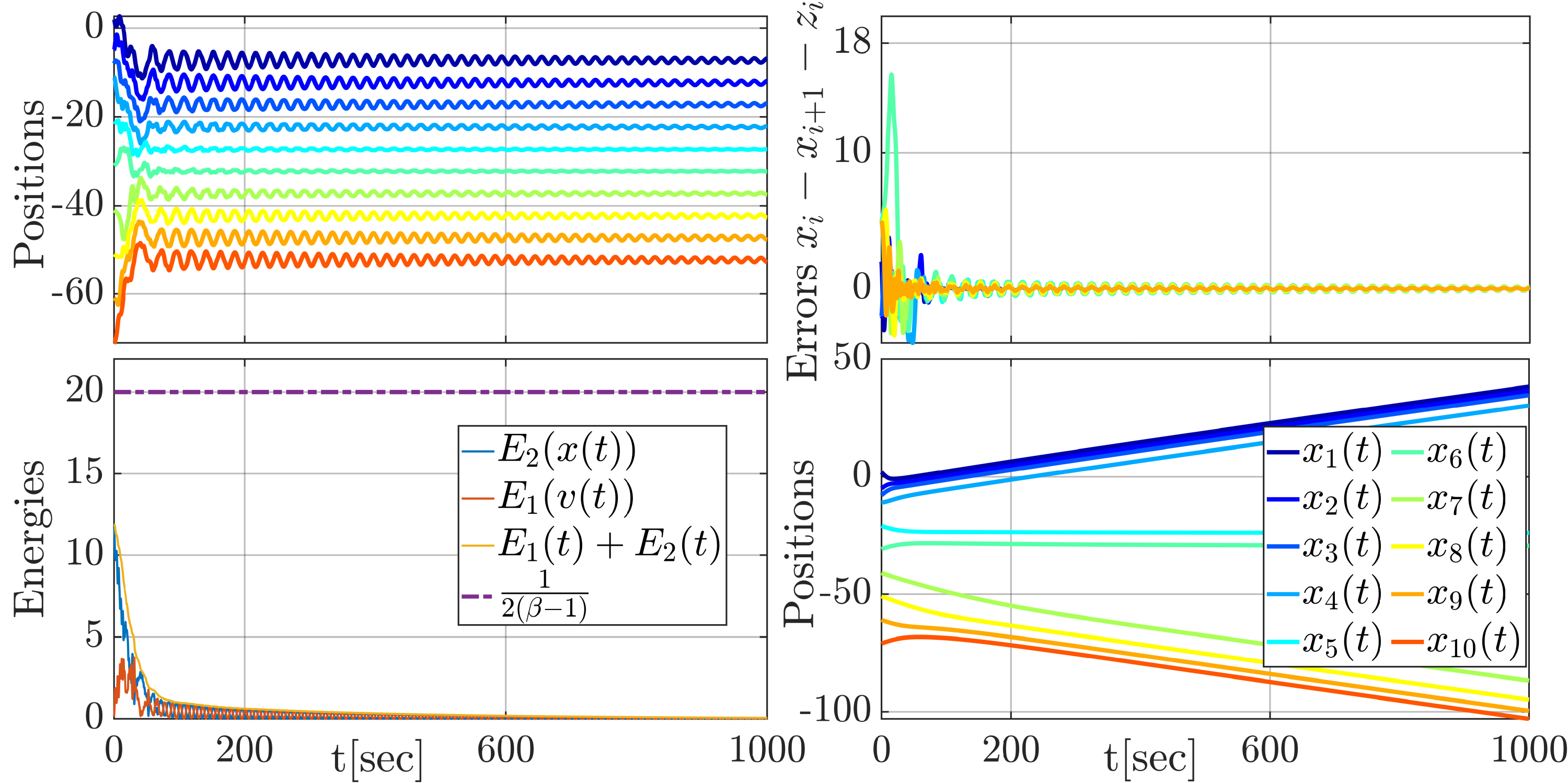}\\
		\caption{10 particles on the line with zero average velocity when $\beta=1.025$.  Left-Top: Positions over time of the particles when the control is used; Left-Bottom: Energy decomposition and flocking condition for $\beta>1$; Right-Top: Errors from the desired formation $x_i-x_{i+1}-z_i$; Right-Bottom: Positions over time of the particles when the control is not used. As predicted by Theorem \ref{thm_1}, given that the initial total energy satisfies \eqref{eq:condT2}, no collisions occur and the desired formation is achieved in steady state.} \label{fig:beta1.5}
	\end{center}
\end{figure}

We now consider $N=10$ agents with $\alpha=2.2$ and $\delta_i=0$ for all $i$. For the same initial conditions with $v_c(0)=0$, we consider two cases: 1) $\beta=4.1$ and 2) $\beta=1.025$, which are presented in Figs. \ref{fig:beta4.1} and \ref{fig:beta1.5} respectively. It can be seen that the parameter $\beta$ influences the value of the initial energies, as noted earlier. However, the kinetic energy is positive and the same in both cases, that is, the agents are not initially moving with the same velocity. For $\beta=4.1$, condition \eqref{main_as} of Theorem \ref{thm_1} is not satisfied, as the value $0.5/(4.1-1)\approx 0.1613$ is less than the initial energy of the system. We can see in Fig. \ref{fig:beta4.1} top-right that the system does not achieve flocking nor the desired formation. Some errors $x_i(t)-x_{i+1}(t)-z_i$ diverge and we can appreciate clustering. On the other hand, for $\beta=1.025$, condition  \eqref{main_as} is met, and we have that the system does achieve the desired formation with all the errors $|x_i(t)-x_{i+1}(t)-z_i|\rightarrow 0$ as the system evolves. Note that the uncontrolled system is plotted at the bottom-right corner in both cases for comparison, that is, with $u_i(t)=0$ for all $t,i$.

\subsection*{Conclusions}\label{concl}
We have presented a control system for platooning composed by a string of agents interacting under nonlinear singular dynamics and a decentralized feedback law. The resulting closed-loop exhibits important features for platooning control, namely, collision-avoidance, velocity flocking, and asymptotic pattern formation. The derivation of rigorous energy estimates allow the characterization of conditions under which the aforementioned features are guaranteed. Energy estimates are governed by: the number of agents in the string, the strength of the control interaction term expressed through the parameter $\beta$ in \eqref{kers}, and the cohesiveness of the initial configuration. In particular, the dependence with respect to the number of agents is a relevant topic of interest for future research. Although our results are asymptotic, we have observed the transient behaviour of the control system and it exhibits similarities to linear-time invariant platooning, namely, slow transients as the number of agents increases. The energy analysis we presented can be extended to study mean field dynamics arising when $N\to\infty$ and the system is characterized by an agent density function \cite{Piccoli2020}. Although the applicability of the mean field framework seems inadequate from a safety viewpoint as collision-avoidance is an eminently microscopic phenomenon, it can be a powerful mathematical method to further understand the large-scale structure of the control system.

\subsection*{Acknowledgements}
YPC has been supported by NRF grant (No. 2017R1C1B2012918) and Yonsei University Research Fund of 2020-22-0505. DK was supported by the UK Engineering and Physical Sciences Research Council (EPSRC) grants and EP/T024429/1. APR was supported by ANID FONDECYT 11221365 grant.

\bibliographystyle{plain}
\bibliography{bibckp}

\begin{thebibliography}{10}

\bibitem{achl}
Shin~Mi Ahn, Heesun Choi, Seung-Yeal Ha, and Ho~Lee.
\newblock On collision-avoiding initial configurations to {C}ucker-{S}male type
  flocking models.
\newblock {\em Commun. Math. Sci.}, 10(2):625--643, 2012.

\bibitem{surveyalbi}
G.~Albi, N.~Bellomo, L.~Fermo, S.-Y. Ha, J.~Kim, L.~Pareschi, D.~Poyato, and
  J.~Soler.
\newblock Vehicular traffic, crowds, and swarms: From kinetic theory and
  multiscale methods to applications and research perspectives.
\newblock {\em Mathematical Models and Methods in Applied Sciences},
  29(10):1901--2005, 2019.

\bibitem{AKKJMLR}
Behzad Azmi, Dante Kalise, and Karl Kunisch.
\newblock Optimal feedback law recovery by gradient-augmented sparse polynomial
  regression.
\newblock {\em Journal of Machine Learning Research}, 22(48):1--32, 2021.

\bibitem{bbck}
Rafael Bailo, Mattia Bongini, Jos{\'e}~A Carrillo, and Dante Kalise.
\newblock Optimal consensus control of the {C}ucker-{S}male model.
\newblock {\em IFAC-PapersOnLine}, 51(13):1--6, 2018.

\bibitem{balch}
Tucker Balch and Ronald~C. Arkin.
\newblock Communication in reactive multiagent robotic systems.
\newblock {\em Autonomous Robots}, 1(1):27--52, 1994.

\bibitem{beaver2021overview}
Logan~E Beaver and Andreas~A Malikopoulos.
\newblock An overview on optimal flocking.
\newblock {\em Annual Reviews in Control}, 51:88--99, 2021.

\bibitem{bf}
Mattia Bongini and Massimo Fornasier.
\newblock Sparse stabilization of dynamical systems driven by attraction and
  avoidance forces.
\newblock {\em Netw. Heterog. Media}, 9(1):1--31, 2014.

\bibitem{bfk}
Mattia Bongini, Massimo Fornasier, and Dante Kalise.
\newblock ({U}n)conditional consensus emergence under perturbed and
  decentralized feedback controls.
\newblock {\em Discrete Contin. Dyn. Syst.}, 35(9):4071--4094, 2015.

\bibitem{ha22mathphys}
Junhyeok Byeon, Seung-Yeal Ha, and Jeongho Kim.
\newblock Asymptotic flocking dynamics of a relativistic cucker--smale flock
  under singular communications.
\newblock {\em Journal of Mathematical Physics}, 63(1):012702, 2022.

\bibitem{byeon2023emergence}
Junhyeok Byeon, Seung-Yeal Ha, and Jeongho Kim.
\newblock Emergence of state-locking for the first-order nonlinear consensus
  model on the real line.
\newblock {\em Kinetic \& Related Models}, 16(3), 2023.

\bibitem{cfpt}
Marco Caponigro, Massimo Fornasier, Benedetto Piccoli, and Emmanuel Tr\'elat.
\newblock Sparse stabilization and control of alignment models.
\newblock {\em Math. Models Methods Appl. Sci.}, 25(3):521--564, 2015.

\bibitem{cch14}
Jos\'e~A. Carrillo, Young-Pil Choi, and Maxime Hauray.
\newblock Local well-posedness of the generalized {C}ucker-{S}male model with
  singular kernels.
\newblock In {\em M{MCS}, {M}athematical modelling of complex systems},
  volume~47 of {\em ESAIM Proc. Surveys}, pages 17--35. EDP Sci., Les Ulis,
  2014.

\bibitem{CCMP}
Jos\'e~A. Carrillo, Young-Pil Choi, Piotr~B. Mucha, and Jan Peszek.
\newblock Sharp conditions to avoid collisions in singular {C}ucker-{S}male
  interactions.
\newblock {\em Nonlinear Anal. Real World Appl.}, 37:317--328, 2017.

\bibitem{cheng2022collision}
Jianfei Cheng, Lining Ru, Xiao Wang, and Yicheng Liu.
\newblock Collision-avoidance, aggregation and velocity-matching in a
  cucker--smale-type model.
\newblock {\em Applied Mathematics Letters}, 123:107611, 2022.

\bibitem{Choi2017}
Young-Pil Choi, Seung-Yeal Ha, and Zhuchun Li.
\newblock Emergent dynamics of the {C}ucker-{S}male flocking model and its
  variants.
\newblock In {\em Active particles. {V}ol. 1. {A}dvances in theory, models, and
  applications}, Model. Simul. Sci. Eng. Technol., pages 299--331.
  Birkh\"auser/Springer, Cham, 2017.

\bibitem{ckpp19}
Young-Pil Choi, Dante Kalise, Jan Peszek, and Andr{\'e}s~A Peters.
\newblock A collisionless singular cucker-{S}male model with decentralized
  formation control.
\newblock {\em SIAM Journal on Applied Dynamical Systems}, 18(4):1954--1981,
  2019.

\bibitem{choi2021one}
Young-Pil Choi and Xiongtao Zhang.
\newblock One dimensional singular cucker--smale model: uniform-in-time
  mean-field limit and contractivity.
\newblock {\em Journal of Differential Equations}, 287:428--459, 2021.

\bibitem{CD}
F.~Cucker and J.~G. Dong.
\newblock Avoiding collisions in flocks.
\newblock {\em IEEE Transactions on Automatic Control}, 55(5):1238--1243, 2010.

\bibitem{cs07}
F.~Cucker and S.~Smale.
\newblock Emergent behavior in flocks.
\newblock {\em IEEE Transactions on Automatic Control}, 52(5):852--862, 2007.

\bibitem{cuckerhuepe}
Felipe Cucker and Cristian Huepe.
\newblock Flocking with informed agents.
\newblock {\em Mathematic{S} in {A}ction}, 1:1--25, 2008.

\bibitem{dalmao}
F.~Dalmao and E.~Mordecki.
\newblock {C}ucker--{S}male flocking under hierarchical leadership and random
  interactions.
\newblock {\em SIAM Journal on Applied Mathematics}, 71(4):1307--1316, 2011.

\bibitem{feng2019}
Shuo Feng, Yi~Zhang, Shengbo~Eben Li, Zhong Cao, Henry~X Liu, and Li~Li.
\newblock String stability for vehicular platoon control: Definitions and
  analysis methods.
\newblock {\em Annual Reviews in Control}, 2019.

\bibitem{HaHaKim}
Seung-Yeal Ha, Taeyoung Ha, and Jong-Ho Kim.
\newblock Asymptotic dynamics for the {C}ucker-{S}male-type model with the
  {R}ayleigh friction.
\newblock {\em J. Phys. A}, 43(31):315201, 19, 2010.

\bibitem{ha20critical}
Seung-Yeal Ha, Zhuchun Li, and Xiongtao Zhang.
\newblock On the critical exponent of the one-dimensional cucker--smale model
  on a general graph.
\newblock {\em Mathematical Models and Methods in Applied Sciences},
  30(09):1653--1703, 2020.

\bibitem{ha2020critical}
Seung-Yeal Ha, Zhuchun Li, and Xiongtao Zhang.
\newblock On the critical exponent of the one-dimensional cucker--smale model
  on a general graph.
\newblock {\em Mathematical Models and Methods in Applied Sciences},
  30(09):1653--1703, 2020.

\bibitem{ha2018first}
Seung-Yeal Ha, Jinyeong Park, and Xiongtao Zhang.
\newblock A first-order reduction of the cucker--smale model on the real line
  and its clustering dynamics.
\newblock {\em Communications in Mathematical Sciences}, 16(7):1907--1931,
  2018.

\bibitem{helbing}
Dirk Helbing and Peter Molnar.
\newblock Social force model for pedestrian dynamics.
\newblock {\em Physical review E}, 51(5):4282, 1995.

\bibitem{hk18}
Michael Herty and Dante Kalise.
\newblock Suboptimal nonlinear feedback control laws for collective dynamics.
\newblock {\em 2018 IEEE 14th International Conference on Control and
  Automation (ICCA)}, pages 556--561, 2018.

\bibitem{huang}
Minyi Huang, Caines PE, and R.~P. Malhame.
\newblock Individual and mass behaviour in large population stochastic wireless
  power control problems: centralized and nash equilibrium solutions.
\newblock {\em 42nd IEEE International Conference on Decision and Control (IEEE
  Cat. No.03CH37475)}, 1:98--103 Vol.1, 2003.

\bibitem{huang2022stochastic}
Qiao Huang and Xiongtao Zhang.
\newblock On the stochastic singular cucker--smale model: Well-posedness,
  collision-avoidance and flocking.
\newblock {\em Mathematical Models and Methods in Applied Sciences},
  32(01):43--99, 2022.

\bibitem{kim2021first}
Jeongho Kim.
\newblock First-order reduction and emergent behavior of the one-dimensional
  kinetic cucker-smale equation.
\newblock {\em Journal of Differential Equations}, 302:496--532, 2021.

\bibitem{lasry}
Jean-Michel Lasry and Pierre-Louis Lions.
\newblock Mean field games.
\newblock {\em Japanese Journal of Mathematics}, 2(1):229--260, 2007.

\bibitem{leslie2024finite}
Trevor~M Leslie and Changhui Tan.
\newblock Finite-and infinite-time cluster formation for alignment dynamics on
  the real line.
\newblock {\em Journal of Evolution Equations}, 24(1):8, 2024.

\bibitem{mp18}
Piotr~B. Mucha and Jan Peszek.
\newblock The {C}ucker--{S}male equation: Singular communication weight,
  measure-valued solutions and weak-atomic uniqueness.
\newblock {\em Archive for Rational Mechanics and Analysis}, 227(1):273--308,
  2018.

\bibitem{olfati}
R.~Olfati-Saber, J.~A. Fax, and R.~M. Murray.
\newblock Consensus and cooperation in networked multi-agent systems.
\newblock {\em Proceedings of the IEEE}, 95(1):215--233, 2007.

\bibitem{pkh10}
Jaemann Park, H~Jin Kim, and Seung-Yeal Ha.
\newblock Cucker-smale flocking with inter-particle bonding forces.
\newblock {\em IEEE transactions on automatic control}, 55(11):2617--2623,
  2010.

\bibitem{perea}
Laura Perea, Gerard G{\'o}mez, and Pedro Elosegui.
\newblock Extension of the cucker-smale control law to space flight formations.
\newblock {\em Journal of guidance, control, and dynamics}, 32(2):527--537,
  2009.

\bibitem{Peters}
Andrés~A. Peters, Richard~H. Middleton, and Oliver Mason.
\newblock Leader tracking in homogeneous vehicle platoons with broadcast
  delays.
\newblock {\em Automatica}, 50(1):64 -- 74, 2014.

\bibitem{piccoli2015control}
Benedetto Piccoli, Francesco Rossi, and Emmanuel Tr{\'e}lat.
\newblock Control of the 1d continuous version of the cucker-smale model.
\newblock In {\em 2015 American Control Conference (ACC)}, pages 1264--1269.
  IEEE, 2015.

\bibitem{Piccoli2020}
Benedetto Piccoli, Andrea Tosin, and Mattia Zanella.
\newblock Model-based assessment of the impact of driver-assist vehicles using
  kinetic theory.
\newblock {\em Zeitschrift f{\"u}r angewandte Mathematik und Physik}, 71:1--25,
  2020.

\bibitem{sepahe04}
Peter Seiler, Aniruddha Pant, and Karl Hedrick.
\newblock Disturbance propagation in vehicle strings.
\newblock {\em IEEE Transactions on automatic control}, 49(10):1835--1842,
  2004.

\bibitem{sumpter}
David~JT Sumpter, Jens Krause, Richard James, Iain~D Couzin, and Ashley~JW
  Ward.
\newblock Consensus decision making by fish.
\newblock {\em Current Biology}, 18(22):1773--1777, 2008.

\bibitem{vicsek}
Tam{\'a}s Vicsek, Andr{\'a}s Czir{\'o}k, Eshel Ben-Jacob, Inon Cohen, and Ofer
  Shochet.
\newblock Novel type of phase transition in a system of self-driven particles.
\newblock {\em Physical review letters}, 75(6):1226, 1995.

\bibitem{wang2019survey}
Ziran Wang, Yougang Bian, Steven~E Shladover, Guoyuan Wu, Shengbo~Eben Li, and
  Matthew~J Barth.
\newblock A survey on cooperative longitudinal motion control of multiple
  connected and automated vehicles.
\newblock {\em IEEE Intelligent Transportation Systems Magazine}, 12(1):4--24,
  2019.

\bibitem{yin2022nonexistence}
Xiuxia Yin, Zhiwei Gao, Zili Chen, and Yichuan Fu.
\newblock Nonexistence of the asymptotic flocking in the cucker$-$smale model
  with short range communication weights.
\newblock {\em IEEE Transactions on Automatic Control}, 67(2):1067--1072, 2022.

\bibitem{yin2020asymptotic}
Xiuxia Yin, Dong Yue, and Zili Chen.
\newblock Asymptotic behavior and collision avoidance in the cucker–smale
  model.
\newblock {\em IEEE Transactions on Automatic Control}, 65(7):3112--3119, 2020.

\bibitem{ZHANG2020201}
Xiongtao Zhang and Tingting Zhu.
\newblock Complete classification of the asymptotical behavior for singular
  {C}-{S} model on the real line.
\newblock {\em Journal of Differential Equations}, 269(1):201--256, 2020.

\bibitem{zhang2023pattern}
Yinglong Zhang.
\newblock Pattern formation in the cucker-smale model.
\newblock {\em Journal of Differential Equations}, 376:204--234, 2023.

\end{thebibliography}
\end{document}